\documentclass[10pt]{article}
\usepackage{amsmath,amssymb,amsfonts}
\usepackage{epsf,color}
\usepackage{fancyhdr}
\usepackage{blindtext}
\newtheorem{theorem}{Theorem}

\newtheorem{lemma}{Lemma}

\date{}

\title{A Note on the Distribution of the Sum of Lengths of the Initial
Longest Increasing Sequences in Cycles of a Random Permutation}
\bigskip

\author{  {\sc Ljuben Mutafchiev}\thanks{E-mail: Ljuben@aubg.edu}\,\,\thanks{Also at: The Institute of Mathematics
and Informatics, Bulgarian Academy of Sciences,
  Sofia 1113, Bulgaria} \\ \small American University in Bulgaria,
  2700 Blagoevgrad, Bulgaria}
  \date{}
\begin{document}
\maketitle

\begin{abstract}
Let $S_n$ be the set of all permutations of $\{1,2,\ldots,n\}$ and
let $\sigma=(\sigma_1,\sigma_2,\ldots,\sigma_n)\in S_n$. The {\it
initial longest increasing sequence} (ILIS) in $\sigma$ has length
$m$ if, for $1\le m\le n-1$, $\sigma_1<\sigma_2<,\ldots,<\sigma_m,
\sigma_m>\sigma_{m+1}$, and has length $n$ if
$\sigma=(1,2,\ldots,n)$. Let $l(\sigma)$ be the length of the ILIS
in $\sigma$. We assume that $\sigma$ is represented in cycle
notation, so that the first number of each cycle is the minimum
number of this cycle. We also assume that $\sigma$ is chosen
uniformly at random from $S_n$, i.e., with probability $1/n!$. Let
$C_n(\sigma)$ be the set of all cycles of $\sigma$. In \cite{M23},
T. Mansour investigated enumerative properties related to lengths of
the ILIS in random permutations represented by the cycle notation.
In particular, he studied the sum of the ILIS' lengths defined by
$s_n=\sum_{c\in C_n(\sigma)} l(c)$ and derived exact and asymptotic
expressions for its expectation and variance. In this note, we
supplement Mansour's results on $s_n$ with a limit theorem. We show
that $s_n$, appropriately normalized, converges weakly to a standard
normal random variable as $n\to\infty$.

\end{abstract}

\vspace{.5cm}



\vspace{.2cm}

\section{Introduction and Statement of the Main Results}

We start with some notation that will be used further. Let $S_n$ be
the set of all permutations of $[n]:=\{1,2,\ldots,n\}$. We introduce
the uniform probability measure $\mathbb{P}$ on the set $S_n$, that
is, we assign the probability $1/n!$ to each $\sigma\in S_n$. The
notation $\mathbb{E}$ stands for the expectation with respect to
$\mathbb{P}$. Clearly, any permutation
\begin{equation}\label{perm}
\sigma=(\sigma_1,\sigma_2,\ldots,\sigma_n)\in S_n
\end{equation}
corresponds to a directed graph $G_\sigma=([n],\{(i,\sigma_i):i\in
[n]\})$, which is a union of disjoint directed cycles. All numerical
characteristics of $\sigma\in S_n$ (and $G_\sigma$) become random
variables (permutation statistics) on the finite probability space
defined above. Some of them, like number of cycles, length of the
longest cycle, number of cycles of length $k$, $1\le k\le n$, do not
depend on the vertex labeling in $G_\sigma$. We are interested here
in permutation statistics, which depend on the labels of the
vertices in $G_\sigma$. A typical statistic of this type is the
length $L_n(\sigma)$ of the longest increasing subsequence in a
random permutation given by (\ref{perm}). It is defined by
\begin{eqnarray} L_n=L_n(\sigma):=\max\{m: \text{there exist}\quad 1\le
i_1<i_2<\ldots<i_m,
\nonumber \\
 \text{such that}\quad
\sigma_{i_1}<\sigma_{i_2}<\ldots<\sigma_{i_m}\}. \nonumber
\end{eqnarray}
Baik, Deift and Johansson \cite{BDJ99} proved a remarkable limit
theorem for $L_n$. They found an appropriate normalization for $L_n$
and showed that it converges weakly to a random variable having the
Tracy-Widom distribution function. It has been recognized for a long
time that, together with the enumerative combinatorics, this problem
also concerns several disparate mathematical areas (e.g., random
matrix theory, representation theory, integrable systems and
statistical mechanics). For more details, we refer the reader to
\cite{AD99, C18, R15}.

In this note, we focus on the length of the {\it initial longest
increasing sequence} (ILIS) of a random permutation, represented as
a union of cycles. This problem was recently studied by Mansour
\cite{M23}. Our work supplements in a probabilistic aspect one of
his main results and may be considered as a continuation of his
study.

To introduce the concept of the ILIS, we consider again the
permutation given by (\ref{perm}). The ILIS in $\sigma$ has length
$m$ if, for $1\le m\le n-1$, $\sigma_1<\sigma_2<,\ldots,<\sigma_m,
\sigma_m>\sigma_{m+1}$, and has length $n$ if
$\sigma=(1,2,\ldots,n)$. For $\sigma\in S_n$, we define
$l=l(\sigma)$ to be the length of the ILIS in $\sigma$. To
illustrate, consider the permutation
\begin{equation}\label{ex}
\tilde{\sigma}=(1,3,5,4,7,6,2)\in S_7.
\end{equation}
Obviously, $l(\tilde{\sigma})=3$ since the ILIS in $\tilde{\sigma}$
is $1,3,5$. For the longest subsequence in $\tilde{\sigma}$, we have
$L_7(\tilde{\sigma})=4$, but note that there are two subsequences of
length $4$ in $\tilde{\sigma}$: $1,3,5,7$ and $1,3,5,6$. Also note
that the longest increasing subsequence in a permutation is in
general not unique.

In \cite{M23}, Mansour showed that if $\sigma$ is written in the
one-line form (\ref{perm}), then the probability mass function and
the expectation of $l(\sigma)$ can be easily computed, using a
direct counting argument. The main results of \cite{M23} are devoted
to permutations that are represented as unions of cycles, where each
cycle is considered as a separate indecomposable (one-cycle)
permutation. To introduce the statistics studied in \cite{M23}, for
$\sigma\in S_n$, we denote by $C_n(\sigma)$ the set of all cycles of
$\sigma$. We suppose that the numbers in each cycle of $C_n(\sigma)$
are written in a linear order, so that its minimum number is the
first number in that order. In addition, if
$C_n(\sigma)=\{c_1,c_2,\ldots,c_d\}$, where $1\le d\le n$, we also
suppose that the cycles $c_j$ are arranged in the following way:
$\min{\{i\in c_1\}}<\min{\{i\in c_2\}}<\ldots<\min{\{i\in c_d\}}$.
Consider now the following two statistics of a random permutation
$\sigma\in S_n$:
$$
\mathcal{L}_n=\mathcal{L}_n(\sigma) :=\max{\{l(c):c\in
C_n(\sigma)\}}
$$
and
\begin{equation}\label{s}
s_n=s_n(\sigma):=\sum_{c\in C_n(\sigma)} l(c).
\end{equation}
The cycle decomposition of the permutation $\tilde{\sigma}$ given by
(\ref{ex}) is $\tilde{\sigma}=(1)(2357)(4)(6)$. Hence
$\mathcal{L}_7(\tilde{\sigma})=4$ and
$s_7(\tilde{\sigma})=1+4+1+1=7$. For fixed $m$, Mansour \cite{M23}
obtained an expression for the exponential generating  function of
the probabilities $\mathbb{P}(\mathcal{L}_n\le m)$, from which he
deduced the exact and asymptotic distributions of $\mathcal{L}_n$.
Mansour has also obtained the exponential generating function of the
expectations $\mathbb{E}(y^{s_n})$ and derived exact and asymptotic
expressions for the expected value and the variance of $s_n$.

Our main goal in this note is to find a proper normalization for
$s_n$, which yields convergence to the standard normal distribution
as $n\to\infty$. In the next section we shall prove the following

\begin{theorem}.
Let
\begin{equation}\label{spr}
s_n^\prime:=\frac{s_n-e\log{n}}{\sqrt{3e\log{n}}}, \quad n\ge 1.
\end{equation}
The random variable $s_n^\prime$ converges in distribution to a
standard normal random variable as $n\to\infty$, i.e.,
\begin{equation}\label{conv}
\lim_{n\to\infty}\mathbb{P}(s_n^\prime\le u)
=\frac{1}{\sqrt{2\pi}}\int_{-\infty}^u e^{-v^2/2} dv, \quad
-\infty<u<\infty.
\end{equation}
\end{theorem}

In the proof of Theorem 1, we apply:
\begin{itemize}
\item{} Mansour's identity for the exponential generating function
of the expectations $\mathbb{E}(y^{s_n})$ (see \cite[Theorem
2.3]{M23});
\item{} Darboux's theorem for asymptotic estimates of the coefficients of
generating functions with algebraic singularities on the circle of
convergence \cite{D78};
\item{} Curtiss' continuity theorem for moment generating functions
\cite{C42}.
\end{itemize}

\section{Proof of Theorem 1}

We start by considering Mansour's identity. We give it as a separate
lemma.

\begin{lemma}
\cite[Theorem 3.2]{M23}. Let
\begin{equation}\label{f}
f_n(y):=\sum_{\sigma\in S_n} y^{s_n(\sigma)}, \quad n\ge 1,\quad
f_0(y):=1,
\end{equation}
and let
$$
H(x,y):=\sum_{n\ge 0}\frac{f_n(y)}{n!}x^n
$$
be the exponential generating function of the functions $f_n(y)$.
Then we have
\begin{equation}\label{id}
H(x,y) =\frac{1}{(1-x)^{1-(1-y)e^y}} e^{(1-y)e^y \int_0^1
\frac{e^{-y(1-x)u}-e^{-yu}}{u}du}.
\end{equation}
\end{lemma}

It is easily seen that Lemma 1 provides a formula for the
exponential generating function of the expectations
$\mathbb{E}(y^{s_n})$. In fact, from the decomposition
$S_n=\cup_{j=1}^n\{\sigma:s_n(\sigma)=j\}$ and (\ref{f}) it follows
that
$$
\frac{f_n(y)}{n!}=\sum_{j=1}^n\mathbb{P}(s_n=j)y^j
=\mathbb{E}(y^{s_n}), \quad n\ge 1.
$$
Hence (\ref{id}) becomes
\begin{equation}\label{idgf}
H(x,y)=\sum_{n\ge 0}\mathbb{E}(y^{s_n})x^n
=\frac{1}{(1-x)^{1-(1-y)e^y}} e^{(1-y)e^y \int_0^1
\frac{e^{-y(1-x)u}-e^{-yu}}{u}du}.
\end{equation}

Further on, we make a convention that $x$ is a variable from the
complex plane and that $y$ assumes real values in a neighborhood of
$y=1$, say,
\begin{equation}\label{del}
1-\delta<y<1+\delta, \quad 0<\delta<e^{-2}.
\end{equation}
The function in the right-hand side of (\ref{idgf}) is of the
following general form:
\begin{equation}\label{algsin}
G(x):=(1-x)^\alpha g(x).
\end{equation}
We shall introduce now the concept of an {\it algebraic singularity}
on the circle of convergence of a function of a complex variable.
Let us assume that $\alpha$ and $g(x)$ from the right-hand side of
(\ref{algsin}) satisfy the following three conditions:

(i) $\alpha\notin\{0,1,\ldots\}$;

(ii) $g(x)$ is a function analytic in a neighborhood of $x=1$;

(iii) $g(1)\neq 0$.

A singularity of a certain function $F(x)$ is called {\it algebraic}
at $x=1$ if $F(x)$ can be written as the sum of a function analytic
in a neighborhood of $1$ and a finite number of terms of the form
(\ref{algsin}). Darboux's theorem \cite{D78} gives asymptotic
expansions of the Taylor coefficients for functions with algebraic
singularities on their circle of convergence. For more details on
singularity analysis of generating functions, we also refer the
reader to \cite[Section 11]{O95} and \cite[Chapter VI]{FS09}. In
(\ref{idgf}), we encounter an algebraic singularity containing only
a single term. Let $x^n[F(x)]$ denote the coefficient of $x^n$ in
the power series expansion of $F(x)$. Darboux's theorem, as given by
Odlyzko in \cite[Theorem 11.7]{O95}, then yields
\begin{equation}\label{dar}
x^n[G(x)]=g(1)\frac{n^{-\alpha-1}}{\Gamma(-\alpha)}
+o(n^{-\alpha-1}),
\end{equation}
where $G(x)$ is given by (\ref{algsin}) and $\Gamma$ denotes Euler's
gamma function. The function from the right-hand side of
(\ref{idgf}), written in terms of notation (\ref{algsin}), has
\begin{equation}\label{ag}
\alpha=(1-y)e^y-1\\ \quad \text{and} \quad g(x)=e^{(1-y)e^y \int_0^1
\frac{e^{-y(1-x)u}-e^{-yu}}{u}du}.
\end{equation}
Here the variable $y$ may be viewed as a parameter satisfying the
inequalities (\ref{del}). It can be easily verified that conditions
(i)--(iii) hold for $\alpha$ and $g(x)$, defined by (\ref{ag}). In
fact, using the first equality of (\ref{ag}) and both inequalities
in (\ref{del}), we see that $\alpha<e^{\delta-1}-1<0$, which proves
condition (i). To prove (ii), we note that in a neighborhood of
$x=1$, the integral
\begin{equation}\label{i}
\int_0^1 \frac{e^{-y(1-x)u}-e^{-yu}}{u}du
\end{equation}
has an expansion as a series of powers of $x-1$ with coefficients,
which are functions of $y$ and are obtained after term-by-term
integration of the power series expansion of the integrand in
(\ref{i}). Therefore $g(x)$, defined by the second equality of
(\ref{ag}), has also such an expansion in powers of $x-1$. This
implies that $g(x)$ is analytic in a neighborhood of $x=1$ and thus
(ii) is proved. The verification of condition (iii) is
straightforward and follows from the exponential representation of
$g(x)$ given by (\ref{ag}).

Our next goal is to apply the general formula (\ref{dar}) for the
$n$th coefficient of $G(x)=H(x,y)$. From (\ref{idgf}) it follows
that
\begin{equation}\label{expect}
x^n[H(x,y)]=\mathbb{E}(y^{s_n}).
\end{equation}
In the right-hand side of (\ref{dar}) we have to deal with $\alpha$
and $g(x)$ defined by (\ref{ag}). To get a convenient expression for
$g(1)$, we set $x=1$ in (\ref{ag}), expand the integrand in a power
series, and then integrate with respect to $u$. We obtain
\begin{equation}\label{gone}
g(1)=e^{(1-y)e^y h(y)},
\end{equation}
where
\begin{equation}\label{h}
h(y)=\sum_{j=1}^\infty (-1)^{j-1}\frac{y^j}{jj!}.
\end{equation}
We are now ready to apply a particular case of Darboux's theorem in
the form given by (\ref{dar}). Combining (\ref{dar}) with
(\ref{ag}), (\ref{expect}) and (\ref{gone}), we observe that
\begin{equation}\label{asyy}
\mathbb{E}(y^{s_n})=\frac{n^{(y-1)e^y}}{\Gamma(1+(y-1)e^y)}
e^{-(y-1)e^y h(y)}(1+o(1)),\quad n\to\infty,
\end{equation}
where $h(y)$ is defined by (\ref{h}) and $o(1)\to 0$ uniformly for
all $y$ satisfying (\ref{del}).

Further on, we shall use the moment generating function $M_n(t)$ of
the random variable $s_n^\prime$, defined by (\ref{spr}). We have
\begin{equation}\label{mgf}
M_n(t):=\mathbb{E}(e^{ts_n^\prime}) =e^{-t\sqrt{(e\log{n})/3}}
\mathbb{E}(e^{ts_n/\sqrt{3e\log{n}}}).
\end{equation}
We intend to apply Curtiss' continuity theorem \cite{C42} to the
sequence $(M_n(t))_{n\ge 1}$ (the statement and proof of this
theorem may be also found in \cite[Chapter I, Section 2]{S97}).
Roughly speaking, Curtiss' theorem asserts that the uniform
convergence in a neighborhood of the origin of the moment generating
functions of a sequence of random variables  yields the convergence
in distribution (weak convergence) of this sequence. Consequently,
we can obtain the convergence in distribution (\ref{conv}) of the
sequence $(s_n^\prime)_{n\ge 1}$ provided we prove that, for a
certain $t_0>0$,
\begin{equation}\label{key}
\lim_{n\to\infty}M_n(t)=e^{t^2/2}, \quad \text{for all $t$ with
$|t|<t_0$},
\end{equation}
where $M_n(t)$ is given by (\ref{mgf}). We also recall that
$e^{t^2/2}$ is the moment generating function of the standard normal
distribution (see, e.g., \cite[Section 5.8]{GS01}).

The proof of (\ref{key}) is based on an asymptotic estimate for the
expectation $\mathbb{E}(e^{ts_n/\sqrt{3e\log{n}}})$ given in the
right-hand side of (\ref{mgf}). Our analysis stems from
(\ref{asyy}), in which we set
\begin{equation}\label{subst}
y=e^{t/\sqrt{3e\log{n}}}.
\end{equation}
After this change of the variable $y$, (\ref{asyy}) remains true
since by (\ref{key}), $|t|$ is bounded and thus, for large $n$ (say,
$n\ge n_0(t_0,\delta))$, the inequalities from (\ref{del}) are again
satisfied with $y$ given by (\ref{subst}). In what follows next,
using (\ref{subst}), we shall deduce suitable asymptotic expansions
for the terms from the right-hand side of (\ref{asyy}) as
$n\to\infty$. We begin by considering $y-1$ and $e^y$. Expanding the
exponent in (\ref{subst}) in powers of $t$, where $t$ is inherited
from (\ref{key}), we get
\begin{equation}\label{yone}
y-1=\frac{t}{\sqrt{3e\log{n}}} +\frac{t^2}{6e\log{n}}
+O(\log^{-3/2}{n}),
\end{equation}

\begin{eqnarray}\label{ey}
& & e^y =\exp{\left(1+\frac{t}{\sqrt{3e\log{n}}}
+\frac{t^2}{6e\log{n}} +O(\log^{-3/2}{n})\right)} \nonumber \\
& & =ee^{t/\sqrt{3e\log{n}}} e^{t^2/6e\log{n}} (1+O(\log^{-3/2}{n}))
\nonumber \\
& & =e\left(1+(1+\frac{t}{\sqrt{3e\log{n}}} +\frac{t^2}{6e\log{n}}
+O(\log^{-3/2}{n})\right) \nonumber \\
& & \times\left(1+\frac{t^2}{6e\log{n}} +O(\log^{-2}{n})\right)
(1+O(\log^{-3/2}{n})) \nonumber \\
& & =e\left(1+\frac{t}{\sqrt{3e\log{n}}} +\frac{t^2}{3e\log{n}}
+O(\log^{-3/2}{n})\right).
\end{eqnarray}
Combining (\ref{yone}) and (\ref{ey}), we obtain
\begin{equation}\label{pr}
(y-1)e^y =e\left(\frac{t}{\sqrt{3e\log{n}}} +\frac{t^2}{2e\log{n}} +
O(\log^{-3/2}{n})\right).
\end{equation}
This formula implies the less precise estimate
$(y-1)e^y=O(\log^{-1/2}{n})$. Hence, by the polynomial approximation
formula for the gamma function \cite[formula 6.1.35]{AS65}, we have
\begin{equation}\label{gam}
\Gamma(1+(y-1)e^y)=1 +O(\log^{-1/2}{n}).
\end{equation}
We return next to the power series $h(y)$, defined by (\ref{h}). It
is clear that $h(y)$ is bounded for all $y$ from a finite interval.
Combining this with (\ref{pr}) gives
$$
(y-1)e^yh(y)=O(\log^{-1/2}{n}),
$$
and thus
\begin{equation}\label{he}
e^{-(y-1)e^yh(y)}=1+O(\log^{-1/2}{n}).
\end{equation}
Finally, note that (\ref{pr}) also implies
\begin{eqnarray}\label{n}
& & n^{(y-1)e^y} =\exp{((y-1)e^y\log{n})} \nonumber \\
& & =\exp{\left(t\sqrt{\frac{e\log{n}}{3}} +\frac{t^2}{2}
+O(\log^{-1/2}{n})\right)}.
\end{eqnarray}
Putting (\ref{gam})--(\ref{n}) into (\ref{asyy}) with $y$ taken from
(\ref{subst}), we get
$$
\mathbb{E}(e^{ts_n/\sqrt{3e\log{n}}})
=(1+o(1))\exp{\left(t\sqrt{\frac{e\log{n}}{3}} +\frac{t^2}{2}
\right)},
$$
which by (\ref{mgf}) implies (\ref{key}). The proof of Theorem 1 is
now completed after an application of Curtiss' continuity theorem.

\end{document}